\documentclass{article}%
\usepackage{boxedminipage}
\usepackage{framed}
\usepackage{color}
\usepackage{fancybox}%
\usepackage{amsmath}%
\usepackage{amsfonts}%
\usepackage{amssymb}%
\usepackage{graphicx}

\makeatother
\definecolor{shadecolor}{cmyk}{0,0,0,0.15}
\definecolor{PGDarkGreen}{rgb}{0,0.5,0}
\begin{document}

\title{Accurate Computation of Laplace Eigenvalues by an Analytical Level Set Method}
\author{Pavel Grinfeld}
\date{Febuary 15th, 2012}
\maketitle

\begin{abstract}
This purpose of this write-up is to share an idea for accurate computation of
Laplace eigenvalues on a broad class of smooth domains. We represent the
eigenfunction $u$ as a linear combination of eigenfunctions corresponding to
the common eigenvalue $\rho^{2}$:%
\[
u\left(  r,\theta\right)  =\sum_{n=0}^{N}P_{n}J_{n}\left(  \rho r\right)  \cos
n\theta,
\]
We adjust the coefficients $P_{n}$ and the parameter $\rho$ so that the zero
level set of $u$ approximates the domain of interest. For some domains, such
as ellipses of modest eccentricity, the coefficients $P_{n}$ decay
exponentially and the proposed method can be used to compute eigenvalues with
arbitrarily high accuracy.

\end{abstract}

\section{Introduction}

The celebrated level set method \cite{OsherLSMBook} is a numerical method for
solving problems with moving interfaces and has been applied to eigenvalue
problems \cite{OsherSantosaLSMOptimization}. The method proposed here also
represents interfaces as level sets of functions. However, our level set
function is a linear combination of global functions, with the coefficients of
the linear combination and possibly additional parameters acting as the
degrees of freedom for deforming the level set. When only a small number of
functions is required for the accurate computation of some quantity, the
proposed method can be viewed as essentially analytical.

The Laplace eigenvalues (with zero or any other type of boundary conditions)
are not available analytically for a $2\times1$ ellipse. (In fact, eigenvalues
for only a handful of shapes have been discovered since Rayleigh's pioneering
"Theory of Sound" \cite{RayleighTheoryOfSound}.) Furthermore, eigenvalues are
difficult to compute with high accuracy. The finite element\ method
\cite{StrangFix} with isoparametric quadratic elements would perhaps be the
most common approach. Combined with Richardson extrapolation, it is an
impressive $5$-th order method. However, in the author's experience, it fails
to deliver the lowest eigenvalue beyond the $13$-th digit. Furthermore, the
finite element method is ineffective in computing the high
(enough)\ eigenvalues since the supporting mesh must resolve the high
frequency oscillations of the corresponding eigenfunctions.

\section{Description of the approach}

We propose a method that computes the lowest eigenvalue on a $2\times1$
ellipse to arbitrary accuracy. The method can be applied to other shapes and
higher eigenvalues. The method has a number of desirable features, including
utmost simplicity, exponential accuracy, ease of analysis and the ability to
produce an exact solution to an approximate problem.

Represent the domain $\Omega$ with boundary $S$ as the interior $\Omega_{u}$
of the zero level set $S_{u}$ of the function $u\left(  r,\theta\right)  $
given by the linear combination
\begin{equation}
u\left(  r,\theta\right)  =\sum_{n=0}^{N}P_{n}J_{n}\left(  \rho r\right)  \cos
n\theta,\label{u}%
\end{equation}
where $J_{n}$ are Bessel functions. The coefficients $P_{n}$ and the parameter
$\rho$ are used as degrees of freedom in fitting the domain $\Omega_{u}$ to
$\Omega$. Since $u\left(  r,\theta\right)  $ given by (\ref{u}) satisfies%
\begin{equation}
-\Delta u=\rho^{2}u,
\end{equation}
then, by definition, $\rho^{2}$ is a Laplace eigenvalue on $\Omega_{u}\ $with
zero boundary conditions and $u\left(  r,\theta\right)  $ is the corresponding
eigenfunction. Note that while $\Omega_{u}$ is an approximation to $\Omega$,
the function $u\left(  r,\theta\right)  $ is an exact eigenfunction on
$\Omega_{u}$.

The objective function measuring the proximity between $S$ and $S_{u}$ can be
defined as%
\begin{equation}
\int_{S}D_{u}^{2}\left(  S\right)  dS,\label{Objective Function}%
\end{equation}
where $D_{u}$ is an appropriately defined distance function that captures the
discrepancy between $S$ and $S_{u}$. The Mathematica code below takes a
pragmatic approach and defines $D\left(  \alpha\right)  $ as%
\begin{equation}
D\left(  \alpha\right)  =r\left(  \alpha\right)  -r_{u}\left(  \alpha\right)
,
\end{equation}
where $r\left(  \alpha\right)  $ is the polar representation of the ellipse
and $r_{u}\left(  \alpha\right)  $ is the polar representation of $S_{u}$. The
minimization can be carried out effectively in a number of ways, including
available generic optimization routines.

As mentioned above, one of the most appealing characteristics of the proposed
method is its utmost simplicity. The Mathematica code that computes the lowest
eigenvalue on the $2\times1$ ellipse is given below in its entirety.

\texttt{1. Goals = \{ AccuracyGoal -%
$>$
75, PrecisionGoal -%
$>$
75, WorkingPrecision -%
$>$
120\};}

\texttt{2. A = 1/2; B = 1;}

\texttt{3. NumberOfTerms = 30;}

\texttt{4. Shape[theta\_] = 1/Sqrt[Cos[theta]\symbol{94}2/A\symbol{94}2 +
Sin[theta]\symbol{94}2/B\symbol{94}2];}

\texttt{5. lsf[rho\_, p\_][r\_, theta\_] := BesselJ[0, rho*r] + Sum[p[[n]]
BesselJ[2 n, rho*r] Cos[2 n theta], \{n, 1, Length[p]\}];}

\texttt{6. r[rho\_, p\_][theta\_] := x /. FindRoot[lsf[rho, p][x, theta] == 0,
\{x, Shape[theta]\}, Evaluate[Goals]];}

\texttt{7. ObjectiveFunction[rho\_, p : \{\_\_?NumericQ\}] := Sum[(r[rho,
p][theta] - Shape[theta])\symbol{94}2, \{theta, 0, 2 Pi - 2 Pi/60, 2
Pi/60\}]/(2 Pi/60) // Sqrt;}

\texttt{8. s = FindMinimum[ObjectiveFunction[rr, pp], \{\{rr, BesselJZero[0,
1]\}, \{pp, ConstantArray[0, NumberOfTerms]\}\}, Evaluate[Goals],
MaxIterations -%
$>$
500, StepMonitor :%
$>$
Print[rr, " ", ObjectiveFunction[rr, pp] // N[\#, 5] \&, " ", Date[]]];}

Code notes:

Line 1 specifies accuracy goals and the working number of digits. Naturally,
higher working precision takes more time. On a circa 2007 desktop, the
presented code takes about $30$ minutes to take the first step and about $10$
minutes for each subsequent step. It takes several steps to reduce the error
by a factor of $10$. These figures are given for computing with $30$ terms in
series (\ref{u}) and $125$ digits of accuracy. The method is substantially
faster when fewer digits are required. For instance, $6$ digits of accuracy
can be obtained with $3$ terms in series (\ref{u}) in a matter of seconds.

Line 2 specifies the semiaxes of the ellipse.

Line 3 specifies the number of terms in equation (\ref{u}). Due to the
symmetry of the ellipse with respect to the $y$-axis, we only use the even
terms in the series (\ref{u}), raising the effective number of terms to $60$.
The variable NumberOfTerms is referenced in Line 8.

Line 4 gives the ellipse in polar coordinates

Line 5 defines the level set function (lsf)\ $u\left(  r,\theta\right)  $.
This line specifies that $u$ also depends on the parameters $\rho$ and $P_{n}%
$. Since the zero level set remains unchanged when the coefficients $P_{n}$
are multiplied by a number, the expansion here assumes that $P_{0}=1$ and
gives $u\left(  r,\theta\right)  $ as%
\begin{equation}
u\left(  r,\theta\right)  =J_{0}\left(  \rho r\right)  +\sum_{n=1}^{N}%
P_{2n}J_{2n}\left(  \rho r\right)  \cos2n\theta.
\end{equation}

Line 6 defines the zero level set in polar coordinates by solving the equation
$u\left(  r,\theta\right)  =0$.

Line 7 evaluates the objective function in equation (\ref{Objective Function}%
). The integral is approximated by a finite sum with $60$ terms. More terms
are needed when greater accuracy is targeted.

Line 8 finds the optimal configuration. The initial configuration is the unit
circle (since all $P_{n}=0$ for $n>1$).

\section{Further thoughts}

\subsection{Hadamard acceleration}

The level set function $u\left(  r,\theta\right)  $ is an exact eigenfunction
on $\Omega_{u}$. This fact can be used increase the accuracy of the eigenvalue
estimate. The Hadamard formula \cite{HadamardMemoire} gives the rate of change
in the eigenvalue $\lambda$ in response to a deformation of the domain
$\Omega$.\ Consider a smooth evolution $S\left(  t\right)  $ parameterized by
$t$. Then the rate of change $d\lambda/dt$ in the eigenvalue is given by%
\begin{equation}
\frac{d\lambda}{dt}=-\int_{S\left(  t\right)  }C\left\vert \nabla
\psi\right\vert ^{2}dS,
\end{equation}
where $C$ is the Hadamard velocity of $S\left(  t\right)  $ and $\psi$ is the
eigenvalue corresponding to $\lambda$.

Since $\rho^{2}$ is an exact eigenvalue on $S_{u}$ and $u\left(
r,\theta\right)  $ is the corresponding exact eigenfunction, the Hadamard
formula can be effectively used by letting $C$ equal the distance between
$S_{u}$ and $S$ along the normal direction to $S_{u}$.

\subsection{Prescribing $C$}

Deforming the surface according to a prescribed $C$ is likely a more effective
way of finding the optimal $u$ than using a generic optimizer. Suppose that
$D_{u}$ is the signed normal displacement between $S$ and $S_{u}$. Then the
differential equation $C=-D_{u}$ on the moving interface $S_{u}$ corresponds
to $S_{u}$ approaching $S$. Therefore, we may vary the parameters $P_{n}$ and
$\rho$ in such a way that the resulting velocity of $S_{u}$ is as close as
possible to $-D_{u}$.

Let us calculate the velocity of the interface $C$ that results when the
parameters of the level set function $u$ are varied in a prescribed way. In a
general setting, suppose that the function $u$ depends on the parameters
$P_{n}$ and we are now including $\rho$ among the $P_{n}$. Refer $u$ arbitrary
curvilinear coordinates $Z_{i}$. Suppose that the parameters $P_{n}$ vary
according to $P_{n}\left(  t\right)  $ and the resulting zero level set has
the equation $Z_{i}\left(  t,S\right)  $, where $S$ are arbitrary coordinates
on the surface $S_{u}$. The functions $Z_{i}\left(  t,S\right)  $ satisfy the
equation%
\begin{equation}
u\left(  P_{n}\left(  t\right)  ,Z_{i}\left(  t,S\right)  \right)  =0
\end{equation}
Differentiating with respect to $t$, we find%
\begin{equation}
\frac{\partial F}{\partial P_{n}}\dot{P}_{n}+\frac{\partial F}{\partial Z_{i}%
}\frac{\partial Z_{i}\left(  t,S\right)  }{\partial t}=0,
\end{equation}
where $\dot{P}_{n}=dP_{n}/dt$ and summations over $n$ and $i$ are implied.
Denote the gradient $\partial F/\partial Z_{i}$ by $\left\vert \nabla
F\right\vert N$ (where $N_{i}$ is the unit normal) and $\partial
Z_{i}/\partial t$ by $V_{i}$:
\begin{equation}
\frac{\partial F}{\partial P_{n}}\dot{P}_{n}+\left\vert \nabla F\right\vert
N_{i}V_{i}=0,
\end{equation}
Since $\sum N_{i}V_{i}=C$, we find
\begin{equation}
C=-\dot{P}_{n}\frac{\partial F}{\partial P^{n}}\left\vert \nabla F\right\vert
^{-1}.\label{C as a series}%
\end{equation}
Therefore, $\dot{P}_{n}$ should be chosen so that the series
(\ref{C as a series}) approximates $C$ as closely as possible -- for example,
in the least squares sense.

\subsection{Enriching the family of functions}

The presented calculation is based on a very limited set of functions
$J_{n}\left(  \rho r\right)  e^{in\theta}$, all centered at the same origin.
Naturally, any Laplace eigenfunction (corresponding to the same eigenvalue)
can be added to the mix. For the problem at hand, it is beneficial to consider
functions of the form $J_{n}\left(  \rho r\right)  e^{in\theta}$, shifted
other poles. Different geometries may utilize other functions.

\bibliographystyle{abbrv}
\bibliography{Classics,ElectronBubbles,FluidFilms,Fluids,InnerCore,Misc,MovingSurfaces,PGrinfeld,Strang,SymbolicCMS,Tensors}

\end{document}